\documentclass[12pt,twoside]{article}
\usepackage{amsfonts}
\usepackage{amssymb}
\usepackage{amsmath}
{\markboth{\sl L. Kamoun and S.
Negzaoui}{\sl On The Harmonic Analysis Associated to the
Bessel-Struve Operator } \thispagestyle{plain}
\newtheorem{prop}{Proposition}[section]
\newtheorem{thm}{Theorem}[section]
\newtheorem{df}{Definition}[section]
\newtheorem{lem}{Lemma}[section]
\newtheorem{rmk}{Remark}[section]
\newtheorem{cor}{Corollary}[section]

\begin{document}
\title{On The Harmonic Analysis Associated to the Bessel-Struve Operator \footnote{\small The authors are supported by the DGRST research project 04/UR/15-02}}
\author{Lotfi Kamoun$^1$ and Selma Negzaoui$^2$\\\\\small $^1$Department of Mathematics,
Faculty of Sciences of Monastir,\\ \small University of Monastir, 5019 Monastir, Tunisia\\ \small
E-mail : kamoun.lotfi@planet.tn \\ \small $^2$ Preparatory
Institute of Engineering Studies of Monastir, \\ \small University of Monastir, 5000 Monastir,
Tunisia
\\ \small E-mail : selma.negzaoui@issatgb.rnu.tn}
\date{}
\maketitle
\begin{abstract} \noindent In this paper, we introduce
the Bessel-Struve transform, we establish an inversion theorem of
the Weyl integral transform associated with this transform, in the
case of half integers, we give a characterization of the range of
$\mathcal{D}(\mathbb{R})$ by Bessel-Struve transform and we prove
a Schwartz-Paley-Wiener theorem on $\mathcal{E}'(\mathbb{R})$.
\end{abstract}
\section{Introduction}
Discrete harmonic analysis associated with Bessel-Struve kernel
was developed by Watson in \cite{Wa}. He gave some results about
"generalised Schl\"omilch series" which represent the Fourier series
associated to Bessel-Struve kernel. Recently, L.Kamoun and M.Sifi
looked to the Bessel-Struve operator
$$\ell_{\alpha}u(x)=\frac{d^2u}{dx^2}(x)+\frac{2\alpha+1}{x}
\left[\frac{du}{dx}(x)- \frac{du}{dx}(0)\right] $$ which has
Bessel-Struve kernel as eigenfunction. In \cite{kam}, they
considered the Intertwining operator associated with Bessel-Struve
operator on $\mathbb{R}$ and expressed its inverse. A.Gasmi and
M.Sifi introduced the Bessel-Struve transform on the dual space of
entire functions on $\mathbb{C}$ in \cite{gs} and in \cite{gs2} in
order to study mean-periodic functions. In this paper, we consider
the Bessel-Struve transform on $L^1_\alpha (\mathbb{R})$ by
$$\mathcal{F}^\alpha_{BS}(f)(\lambda)=\int_{\mathbb{R}}f(x)\,S_{-i\lambda}^\alpha(x)\,
|x|^{2\alpha+1}\,dx$$ where the Bessel-Struve kernel
$S_\lambda^\alpha$ is given by $$\forall x\in
\mathbb{R},\;S^\alpha_\lambda(x)=j_{\alpha}(i\lambda
x)-ih_{\alpha}(i\lambda x) $$ If $f$ is an even function, this
operator coincides with the Hankel transform defined by $$
\mathcal{H}_\alpha(f)(x)=\int_0^{+\infty} f(t)\,
j_\alpha(xt)\,t^{2\alpha +1}\, dt$$ The outline of the content of
this paper is as follows

 In section 2: We introduce the
Bessel-Struve transform on $L^1_\alpha (\mathbb{R})$, we give a
necessary condition for a function to be the range of
$\mathcal{F}^\alpha_{BS}$ of a function in
$L^1_\alpha(\mathbb{R})$. Besides, we prove that this operator is
symmetric.

In section 3: we deal with the Weyl integral transform associated
to Bessel-Struve operator.\\ In the beginning, we define the dual
operator $\chi_\alpha^\ast$ of the generalized Riemann Liouville
integral $\chi_\alpha$ introduced by L.Kamoun and M.Sifi in
\cite{kam}. We use this operator to introduce Weyl integral
associated to Bessel-Struve operator that we denote $W_\alpha$ by
$$ \forall\; f\in \mathcal{D}(\mathbb{R})\,,\qquad\qquad
\chi^\ast_\alpha\,T_{A\,f}=T_{W_\alpha(f)}$$ where $T_f$ designates
the distribution defined by the function $f$.\\ The Weyl transform
verifies the following relation : $$ \forall\, f\in
\mathcal{D}(\mathbb{R})\,,\qquad
\mathcal{F}^\alpha_{BS}(f)=\mathcal{F}\circ W_\alpha(f)$$ where
$\mathcal{F}$ is the classical Fourier transform.\\ Then, we
establish an inversion theorem of Weyl integral on
$\mathcal{D}(\mathbb{R})$ and we characterize the range of
$\mathcal{D}(\mathbb{R})$ by $W_\alpha$ that we denote
$\Delta_\alpha(\mathbb{R})$. Furthermore we give the expression of
its inverse denoted $\displaystyle
V_{\alpha|\Delta_\alpha(\mathbb{R})}$. Weyl integral associated to
Bessel-Struve operator doesn't have the same properties like in the
classical harmonic analysis: It presents a difficulty on a
singularity in 0 and it doesn't save the space
$\mathcal{D}(\mathbb{R})$.

In section 4: we give some properties of
$\mathcal{F}^\alpha_{BS}(\mathcal{D}(\mathbb{R}))$ for all
$\displaystyle \alpha>-\frac{1}{2}$, we prove a Paley-Wiener
theorem of Bessel-Struve transform in the case $\displaystyle
\alpha=\frac{1}{2}$. Then, by induction, we deduce the result for
half integers. Finally, we give the Paley-Wiener theorem for the
space of distributions with bounded support called
Schawartz-Paley-Wiener theorem.\\ Throughout the paper, we denote
:
\begin{itemize}
  \item $\mathbb{R}^\ast=\mathbb{R}\backslash\{0\}$
  \item $\mathbb{N}^\ast=\mathbb{N}\backslash\{0\}$
\end{itemize}

Finally, we do not forget to specify that the Bessel-Struve operator
and Bessel-Struve kernel were introduced by K. Trim\`eche in an
unpublished paper\cite{t4}.
\section{Bessel-Struve transform}
In this section, we define the Bessel-Struve transform and we give
some results similar to those in the classical harmonic analysis.

First, we consider the differential operator
$\displaystyle\ell_{\alpha},\;\alpha>-{\frac{1}{2}}$, defined on
$\mathbb{R}$ by
\begin{equation}
\displaystyle\ell_{\alpha}u(x)=\frac{d^2u}{dx^2}(x)+\frac{2\alpha+1}{x}
\left[\frac{du}{dx}(x)- \frac{du}{dx}(0)\right]
\end{equation}
with an infinitely differentiable function $u$ on $\mathbb{R}$.
This operator is called Bessel-Struve operator.

For $\lambda\in\mathbb{C}$, the differential equation : $$\left\{
\begin{array}{ll}\displaystyle{\ell_\alpha u(x) =
\lambda^{2}u(x)}\, &\\ \displaystyle{u(0) = 1\;,\;u'(0)=\frac{
\lambda \Gamma(\alpha+1)}{\sqrt{\pi}\,\Gamma(\alpha+3/2)}}\,.&
\end{array}\right.$$
possesses a unique solution denoted $S_\lambda^\alpha$. This
eigenfunction, called the Bessel-Struve kernel, is given by :
\begin{equation}\label{sjh}
S^\alpha_\lambda(x)=j_{\alpha}(i\lambda x)-ih_{\alpha}(i\lambda x)
\end{equation}
 where $j_\alpha$ and $h_\alpha$ are
respectively the normalized Bessel and Struve functions of index
$\alpha$. Those kernels are given as follows :
$$j_\alpha(z)=2^\alpha
\Gamma(\alpha+1)z^{-\alpha}J_\alpha(z)=\Gamma(\alpha+1)\sum_{n=0}^{+\infty}\frac{(-1)^n(z/2)^{2n}}{n!\Gamma(n+\alpha+1)}$$
and $$h_\alpha(z)=2^\alpha
\Gamma(\alpha+1)z^{-\alpha}\mathbf{H}_\alpha(z)=\Gamma(\alpha+1)\sum_{n=0}^{+\infty}
\frac{(-1)^n(z/2)^{2n+1}}{\Gamma(n+\frac{3}{2})\Gamma(n+\alpha+\frac{3}{2})}$$
(for more details one can see \cite{Wa})\\ The kernel
$S_\lambda^\alpha$ possesses the following integral representation
:
\begin{equation}\label{bsk}
\forall x\in \mathbb{R},\;\forall \lambda \in \mathbb{C},\quad
S^\alpha_\lambda(x)=\frac{2\Gamma(\alpha+1)}{\sqrt{\pi}\
\Gamma(\alpha+\frac{1}{2})}\int^{1}_{0}(1-t^{2})^{\alpha-
\frac{1}{2}}e^{\lambda xt}dt
\end{equation}

 We denote by $L^1_\alpha(\mathbb{R})\,$, the space of
real-valued functions $f$, measurable on $\mathbb{R}$ such that
$$\|f\|_{1\,,\,\alpha}=\int_{\mathbb{R}}|f(x)|\,d\mu_\alpha(x)
<+\infty\,,$$ where $$d\mu_\alpha(x)=A(x)\,dx\;\quad
\mathrm{and}\;\quad A(x)=|\,x|^{2\alpha+1}$$
\begin{df}\label{fbs}
We define the Bessel-Struve transform on $L^1_\alpha(\mathbb{R})$
by
\begin{equation}\label{f} \forall \,\lambda\in \mathbb{R},\qquad
\mathcal{F}_{B,S}^\alpha(f)(\lambda)=\int_{\mathbb{R}}f(x)\,S_{-i\lambda}^\alpha(x)\,d\mu_\alpha(x)
\end{equation}
\end{df}
\begin{prop}\label{salpha} {\rm The kernel $S_\lambda^\alpha$ has a unique extension to
$\mathbb{C}\times\mathbb{C}$. It satisfies the following
properties :}
\begin{description}
  \item[(i)] $\forall\lambda\in\mathbb{C},\ \ \forall z\in\mathbb{C},\qquad\qquad
S_{-i\lambda}^\alpha(z)=S_{-iz}^\alpha(\lambda)$
  \item[(ii)] $\forall\lambda\in\mathbb{C},\ \ \forall
z\in\mathbb{C},\qquad\qquad
S_{-\lambda}^\alpha(z)=S_\lambda^\alpha(-z)$
  \item[(iii)] $\displaystyle \forall n\in \mathbb{N},\;\forall\lambda\in
\mathbb{R},\; \forall x\in \mathbb{R},\quad
\left|\frac{d^n}{dx^n}S_{i\lambda}^\alpha(x)\right|\leq
|\lambda|^n$\\ \rm{ In particular, we have for all $\lambda$ and
$x$ in $\mathbb{R}$}, $\quad|S_{i\lambda}^\alpha(x)|\leq 1$
  \item[(iv)] For all $x\in \mathbb{R}^\ast$,\;
  $\displaystyle\lim_{\lambda\rightarrow+\infty}S_{-i\lambda}^\alpha(x)=0$
\end{description}
\end{prop}
{\bf Proof.} The relation (\ref{sjh}) yields directly (i) and
(ii). \\ (iii) For $n\in \mathbb{N}$, we have $$
\frac{d^n}{dx^n}S_{i\lambda}^\alpha(x)
=\frac{2\Gamma(\alpha+1)}{\sqrt{\pi}\Gamma(\alpha+\frac{1}{2})}(i\lambda)^n\int_0^1(1-t^2)^{\alpha-\frac{1}{2}}t^ne^{i\lambda
xt}dt$$ So we find
$\displaystyle|\frac{d^n}{dx^n}S_{i\lambda}^\alpha(x)
|\leq\frac{2\Gamma(\alpha+1)}{\sqrt{\pi}\Gamma(\alpha+\frac{1}{2})}|\lambda|^n\int_0^1(1-t^2)^{\alpha-\frac{1}{2}}dt$
\\ (iv) From the asymptotic expansion of $J_\alpha$ and $\mathbf{H}_\alpha$ p.199 and p.333 in \cite{Wa}, we deduce that $\displaystyle\lim_{\lambda\rightarrow+\infty}j_\alpha(\lambda
  x)=0$ and $\displaystyle\lim_{\lambda\rightarrow+\infty}h_\alpha(\lambda
  x)=0$. Then by relation (\ref{sjh})
$\displaystyle\lim_{\lambda\rightarrow+\infty}S_\alpha(-i\lambda
  x)=0$ \hfill$\square$
\begin{thm}\label{fbsl1}
Let $f$ be a function in $L^1_\alpha(\mathbb{R})$ then
$\mathcal{F}_{B,S}^\alpha(f)$ belongs to $C_0(\mathbb{R})$,\\
where $C_0(\mathbb{R})$ is the space of continuous functions
having 0 as limit in the infinity. Furthermore,
\begin{equation}
\|\mathcal{F}_{B,S}^\alpha(f)\|_\infty\leq\,\|f\|_{1,\alpha}
\end{equation}
\end{thm}
{\bf Proof.}  It's clear that $\mathcal{F}_{B,S}^\alpha(f)$ is a
continuous function on $\mathbb{R}$.\\ From proposition
\ref{salpha}, we get for all $\displaystyle x\in \mathbb{R}^\ast$,
$\displaystyle\lim_{\lambda\rightarrow+\infty}f(x)S_{-i\lambda}^\alpha(x)=0\quad$
and \\ $\quad\displaystyle |f(x)S_{-i\lambda}^\alpha(x)\,A(x)|\leq
|f(x)\,A(x)|$\\ Since $f$ be in $L^1_\alpha(\mathbb{R})$, we
conclude, using the dominated convergence theorem, that
$\mathcal{F}_{B,S}^\alpha(f)$ belongs to $C_0(\mathbb{R})$ and
$$\|\mathcal{F}_{B,S}^\alpha(f)\|_\infty\leq\,\|f\|_{1,\alpha}$$
\hfill$\square$
 \begin{thm} Let $f$ and $g$ in $L^1_\alpha(\mathbb{R})$, we have
\begin{equation}\label{lt}
\int_{\mathbb{R}}\mathcal{F}_{B,S}^\alpha(f)(x)\,g(x)\,d\mu_\alpha(x)=
\int_{\mathbb{R}}\mathcal{F}_{B,S}^\alpha(g)(x)\,f(x)\,d\mu_\alpha(x)
\end{equation}
 \end{thm}
{\bf Proof.} We take
$\phi(t,x)=f(t)g(x)S_{-ix}^\alpha(t)A(t)A(x)$.
\\ Since $f$ and $g$ are both in $L^1_\alpha(\mathbb{R})$ and
according to (iii) of proposition \ref{salpha} and
Fubini-Tonelli's theorem, we deduce that
$\displaystyle\int_{\mathbb{R}^2}|\phi(t,x)|\,dt\,dx\leq\|f\|_{\alpha,1}\|g\|_{\alpha,1}$.
\\ We obtain relation (\ref{lt}) by Fubini's theorem.
 \hfill$\square$
\section{Weyl integral transform}
\subsection{Bessel-Struve intertwining operator and its dual}
$\mathcal{E}(\mathbb{R})$ designates the space of infinitely
differentiable functions on $\mathbb{R}$\,.\\ The Bessel-Struve
intertwining operator on $\mathbb{R}$ denoted $\chi_\alpha$,
introduced by L.Kamoun and M.Sifi in \cite{kam} is defined by:
\begin{equation}
\chi_\alpha(f)(x)=\frac{2\Gamma(\alpha+1)}{\sqrt{\pi}\Gamma(\alpha+\frac{1}{2})}
\int_0^1(1-t^2)^{\alpha-\frac{1}{2}}f(xt)\,dt\;\,,\;\;\;
f\in\mathcal{E}( \mathbb{R})
\end{equation}
\begin{rmk}\begin{equation}\label{schi}
\forall\,x\in \mathbb{R}\,,\;\forall\lambda\in \mathbb{C}\,,\qquad
S_\lambda^\alpha(x)\,=\,\chi_\alpha (e^{\lambda \,.})(x)
\end{equation}
\end{rmk}
\begin{df}
we define the operator $\chi_\alpha^\ast$ on
$\mathcal{E}'(\mathbb{R})$ by
\begin{equation}\label{def1}
<\chi_\alpha^\ast(T),f>=<T,\chi_\alpha f>\qquad,\;\;
f\in\mathcal{E}( \mathbb{R})
\end{equation}
\end{df}
 We denote
$\displaystyle \frac{d}{dx^2}=\frac{1}{2x}\frac{d}{dx}$. The
following theorem is proved in \cite{kam}
\begin{thm}\label{qui-1}
The operator $\chi_\alpha$, $\alpha>-\frac{1}{2}$, is a
topological isomorphism from $\mathcal{E}( \mathbb{R})$ onto
itself. The inverse operator $\chi_\alpha^{-1}$ is given for all
$f\in \mathcal{E}( \mathbb{R})$ by
\begin{description}
  \item[(i)] if $\alpha=r+k\,,\, k\in \mathbb{N}\,,\, -\frac{1}{2}<r<\frac{1}{2} $
$$\hspace*{-1cm}\chi_\alpha^{-1}f(x)=\frac{2\sqrt{\pi}}{\Gamma(\alpha+1)\Gamma(\frac{1}{2}-r)}
x(\frac{d}{dx^2})^{k+1} \left[\int_0^x
(x^2-t^2)^{-r-\frac{1}{2}}|t|^{2\alpha+1}f(t)\,dt\right]$$
  \item[(ii)] if $\alpha=\frac{1}{2}+k\,,\, k\in \mathbb{N}$
$$\chi_\alpha^{-1}f(x)=\frac{2^{2k+1}k!}{(2k+1)!}\,\,x\,(\frac{d}{dx^2})^{k+1}
  \big{(}x^{2k+1}f(x)\,\big{)},\;\;x\in\mathbb{R}  $$
\end{description}
\end{thm}
\begin{cor}\label{tqui-1}
$\chi_\alpha^\ast$ is an isomorphism from
$\mathcal{E}'(\mathbb{R})$ into itself.
\end{cor}
{\bf Proof.} Since $\chi_\alpha$ is an isomorphism from
$\mathcal{E}( \mathbb{R})$ into itself, we deduce the result by
duality. \hfill$\square$
\begin{prop}
For $f\in\mathcal{D}(\mathbb{R})$, the distribution
$\chi_\alpha^\ast  T_{Af}$ is defined by the function $W_\alpha f$
having the following expression
\begin{equation}\label{w}
 W_\alpha f (y)=\frac{2\Gamma(\alpha+1)}{\sqrt{\pi}\,\Gamma(\alpha+\frac{1}{2})}
 \int_{|y|}^{+\infty}(x^2-y^2)^{\alpha-\frac{1}{2}}\,
 x\,f(sgn(y)x)\,dx\;,\;\; y\in \mathbb{R}^\star
\end{equation}
called Weyl integral associated to Bessel-Struve operator.
\end{prop}
{\bf Proof. } Let $g \in \mathcal{E}(\mathbb{R})$,\\
$\displaystyle <\chi_\alpha^\star\, T_{Af},g>\,=\, <
T_{Af},\chi_\alpha g>$\\
\hspace*{2.5cm}$\displaystyle=\frac{2\Gamma(\alpha+1)}{\sqrt{\pi}
\Gamma(\alpha+\frac{1}{2})}\int_0^{+\infty}(\int_0^x(x^2-t^2)^{\alpha-\frac{1}{2}}g(t)\,dt)x\,f(x)\,dx$\\
\hspace*{2.5 cm} $ -\displaystyle
\frac{2\Gamma(\alpha+1)}{\sqrt{\pi} \Gamma(\alpha+\frac{1}{2})}
\int_{-\infty}^0(\int_x^0(x^2-t^2)^{\alpha-\frac{1}{2}}g(t)\,dt)x\,f(x)\,dx$
\\ Using Fubini's theorem, a change of variable and Chasles relation
, we obtain for all $g \in \mathcal{E}(\mathbb{R})$
$$\displaystyle <\chi_\alpha^\star\,
T_{Af},g>\,=\,\int_{\mathbb{R}}\left(\frac{2\Gamma(\alpha+1)}{\sqrt{\pi}
\Gamma(\alpha+\frac{1}{2})}\int_{|t|}^{+\infty}(x^2-t^2)^{\alpha-\frac{1}{2}}x\,f(sgn(t)x)\,dx\right)g(t)\,dt
$$ Therefore, the distribution $\chi_\alpha^\star\, T_{Af}$ is
represented by the function $W_\alpha f$
 \hfill$\square$
\begin{rmk}
Let $f\in \mathcal{E}(\mathbb{R})$ and  $ g\in
\mathcal{D}(\mathbb{R})$. The operator $\chi_\alpha$ and
$W_\alpha$ are related by the following relation
\begin{equation}\label{wqui}
  \int_{\mathbb{R}}\chi_\alpha f(x)\, g(x) A(x)\, dx\,=\,  \int_{\mathbb{R}} f(x)\,W_\alpha g(x) dx
\end{equation}

\end{rmk}
\begin{lem}\label{suppr}
Let $a>0$ and $f\in \mathcal{D}(\mathbb{R})$ with support included
in $[-a,a]$. Then $W_\alpha f$ is infinitely differentiable on
$\mathbb{R}^\star$ and $supp\,(W_\alpha f)$ is included in
$[-a,a]$.\\ Furthermore,  for all $x\in \mathbb{R}^\star$ and
$n\in \mathbb{N},$
\begin{equation}\label{dw}
(W_\alpha f)^{(n)}(x)=\sum_{k=0}^n
\frac{c_\alpha\,(sgn(x))^k}{x^n}
\int_{|x|}^{+\infty}(y^2-x^2)^{\alpha-\frac{1}{2}}y^{k+1}f^{(k)}(y
\, sgn(x))\, dy
\end{equation} where $$
c_\alpha=\frac{2\Gamma(\alpha+1)\,C_n^k\Gamma(2\alpha+2)}{\sqrt{\pi}\,
\Gamma(\alpha+\frac{1}{2})\,\Gamma(2\alpha+2-n+k)}.$$
\end{lem}
{\bf Proof.} Let $f\in\mathcal{D}(\mathbb{R})$ such that
 $supp\,(f)\subseteq [-a,a]$.\\ By change of variable $W_\alpha f$ can be written
\begin{equation}\label{w2}
W_\alpha f(x)=
\frac{2\Gamma(\alpha+1)}{\sqrt{\pi}\,\Gamma(\alpha+\frac{1}{2})}
|x|^{2\alpha+1}\, \int_1^{+\infty}(t^2-1)^{\alpha-\frac{1}{2}}\,
 t\,f(tx)\,dt\quad,\quad x\in \mathbb{R}^\star
\end{equation}
Then, $W_\alpha
 f(x)=A(x)\,\psi(x)$ where
 $$\displaystyle
 \psi(x)=\frac{2\Gamma(\alpha+1)}{\sqrt{\pi}\,\Gamma(\alpha+\frac{1}{2})}
 \int_1^{\frac{a}{|x|}}(t^2-1)^{\alpha-\frac{1}{2}}\,
 t\,f(tx)\,dt\;\;,\quad x\in \mathbb{R}^\ast$$
 If $|x|>a$ then $|tx|>a$ and $f(tx)=0$ so $\psi(x)=0$.
 Consequently, we have  $supp\,(\psi)\subseteq [-a,a]$.
\\ Since $A$ and $\psi$ are both infinitely differentiable on
$\mathbb{R}^\ast$, we have $W_\alpha f $ is infinitely
differentiable on $\mathbb{R}^\ast$ and for all $\displaystyle
n\in \mathbb{N},$ $$(W_\alpha f)^{(n)}(x)=\sum_{k=0}^nC_n^k \,
A^{(n-k)}(x)\, \psi^{(k)}(x)$$ By a change of variable one obtains
$$\displaystyle
\psi^{(k)}(x)=|x|^{-2\alpha-k-1}\frac{2\Gamma(\alpha+1)}{\sqrt{\pi}\Gamma(\alpha+\frac{1}{2})}\int_{|x|}^{a}
(y^2-x^2)^{\alpha-\frac{1}{2}}\, y^{k+1}\,f^{(k)}(sgn(x)y)\,dy $$
Therefore \\ $\displaystyle(W_\alpha f)^{(n)}(x)=$
$$\frac{a_\alpha}{x^n}\sum_{k=0}^n\frac{C_n^k(sgn(x))^k\Gamma(2\alpha+2)}{\Gamma(2\alpha+2-n+k)}
\int_{|x|}^{a}(y^2-x^2)^{\alpha-\frac{1}{2}}y^{k+1}f^{(k)}(y\,sgn(x))\,
dy $$ where \quad $\displaystyle
a_\alpha=\frac{2\Gamma(\alpha+1)}{\sqrt{\pi}\,
\Gamma(\alpha+\frac{1}{2})}$ \hfill$\square$
\begin{cor}\label{xnw} Let $f$ be in $\mathcal{D}(\mathbb{R})$, we have
for all $n\in \mathbb{N}$, $y^n(W_\alpha f)^{(n)}$ belongs to
$L^1(\mathbb{R})$.
\end{cor}
{\bf Proof.} From relation (\ref{dw}), we deduce that
$$\displaystyle \lim_{\substack{ y\rightarrow 0 \\ y>0 }}\;y^n
(W_\alpha f)^{(n)}(y)\quad \mathrm{and} \quad \lim_{\substack{
y\rightarrow 0
\\ y<0 }}\; y^n (W_\alpha f)^{(n)}(y)$$ exist.\\ Consequently we obtain
for all $n\in \mathbb{N}$, $y^n(W_\alpha f)^{(n)}$ belongs to
$L^1(\mathbb{R})$.\hfill$\square$
\begin{prop}\label{wl1}
$W_\alpha $ is a bounded operator from $L^1_\alpha(\mathbb{R})$ to
$L^1(\mathbb{R})$
\end{prop}
{\bf Proof.} Let $f\in L^1_\alpha(\mathbb{R})$ then $\displaystyle
\int_{\mathbb{R}}|f(x)|\,|x|^{2\alpha+1}\,dx<+\infty$
\\  Using Fubini-Tonelli's theorem and a change of variable, we get : \\
$\displaystyle\int_{\mathbb{R}}\left(\int_1^{+\infty}(t^2-1)^{\alpha-\frac{1}{2}}t|f(ty)|\,dt\right)|y|^{2\alpha+1}\,dy
\, = $  \\ \hspace*{3cm}
$\displaystyle\left(\int_{\mathbb{R}}|f(x)||x|^{2\alpha+1}\,dx\right)
\left(\int_1^{+\infty}(t^2-1)^{\alpha-\frac{1}{2}}t^{-2\alpha-1}\,
dt \right)<+\infty  $ \\ Invoking relation(\ref{w2}), we deduce
that $\displaystyle\int_\mathbb{R}|W_\alpha
f(y)|\,dy\quad<\,+\infty\;\;$ \\ and $\; \|W_\alpha(f)\|_1\;\leq
C\; \|f\|_{1,\alpha}$
     \hfill$\square$
\begin{prop}\label{fw} We have
\begin{equation}\label{wf}
\forall f\in \mathcal{D}(\mathbb{R})
\,,\qquad\mathcal{F}_{B,S}^\alpha(f)\,=\;\mathcal{F}\circ
W_\alpha(f)
\end{equation}
 where $\mathcal{F}$ is the classical Fourier
transform defined on $L^1(\mathbb{R})$ by $$\displaystyle
\mathcal{F}(g)(\lambda)=\int_\mathbb{R}g(x)e^{-i\lambda x}dx$$
\end{prop}
{\bf Proof.} Let $f\in \mathcal{D}(\mathbb{R})$. Using
relation(\ref{wqui}) and relation (\ref{schi}), we obtain
$$\mathcal{F}\circ
W_\alpha(f)(x)=\int_\mathbb{R}f(t)\,\chi_\alpha(e^{-ixt})\,A(t)\,dt=\mathcal{F}_{B,S}^\alpha(f)(x)
$$ \hfill$\square$
\subsection{Inversion of Weyl integral}
\begin{lem}
Let $g\in \mathcal{E}( \mathbb{R}^\ast) $, $m$ and $p$ are two
integers nonnegative\,,
 we have
\begin{equation}\label{a}
\forall \,x\in \mathbb{R}^\ast,\qquad\quad
\left(\frac{d}{dx^2}\right)^p(x^{m}\,g(x))=\sum_{i=0}^p\beta_i^p
\, x^{m-2p+i}\,g^{(i)}(x)
\end{equation}
where $\beta_i^p$ are constants depending on $i$, $p$ and $m$.
\end{lem}
{\bf Proof.} We will proceed by induction. The relation (\ref{a})
is true for $p=0$.
\\ Suppose that (\ref{a}) is true at the order
$p\geq 0$ then\\
$\displaystyle(\frac{d}{dx^2})^{p+1}(x^m\,g(x))=\frac{d}{dx^2}(\sum_{i=0}^p\beta_i^p
\, x^{m-2p+i}\,g^{(i)}(x))$\\ \hspace*{3cm} $\displaystyle
=\frac{1}{2x}\sum_{i=0}^p \beta_i^p \,
[(m+i-2p)x^{m-2p+i-1}\,g^{(i)}(x)+x^{m+i-2p}\,g^{(i+1)}(x)]$
\\ \hspace*{3cm} $\displaystyle
=\sum_{i=1}^p\frac{1}{2}[(m+i-2p) \beta_i^p
+\frac{1}{2}\beta_{i-1}^p] x^{m+i-2(p+1)}\,g^{(i)}(x)\,+$\\
\hspace*{3,5cm}$\displaystyle\frac{1}{2}
\beta_0^p(m-2p)x^{m-2(p+1)}g(x)+\frac{1}{2}\beta_p^p\,
x^{m-(p+1)}g^{(p+1)}(x)$
\\ \hspace*{3cm} $\displaystyle
=\sum_{i=0}^{p+1}\beta_i^{p+1}x^{m+i-2(p+1)}\,g^{(i)}(x)$\\ where
$$\beta_{p+1}^{p+1}=\frac{1}{2}\beta_p^p\,,\;
\beta_0^{p+1}=\frac{1}{2}\beta_0^p(m-2p)$$ and $$\forall 1\leq
i\leq p
\,,\;\beta_i^{p+1}=\frac{1}{2}(m+i-2p)\beta^p_i+\frac{1}{2}\beta^p_{i-1}$$
\hfill$\square$

We designate by $\mathcal{K}_0$ the space of functions $f$
infinitely differentiable on $\mathbb{R}^\ast$ with bounded
support and verifying for all $n\in \mathbb{N}$,$$\displaystyle
\lim_{\substack{ y\rightarrow 0 \\ y>0 }}\;y^n (W_\alpha
f)^{(n)}(y)\quad \mathrm{and} \quad \lim_{\substack{ y\rightarrow
0
\\ y<0 }}\; y^n (W_\alpha f)^{(n)}(y)$$ exist.
\begin{prop}\label{v}
Let $f$ a function in $\mathcal{K}_0$. Then the distribution
$(\chi_\alpha^\ast)^{-1}T_f$ is defined by the function denoted
$A\, V_\alpha f$, where $V_\alpha f$ has the following expression
\begin{description}
  \item[(i)]If $\alpha=k+\frac{1}{2},\, k\in \mathbb{N}$
  $$V_\alpha f
  (x)=(-1)^{k+1}\frac{2^{2k+1}k!}{(2k+1)!}\left(\frac{d}{dx^2}\right)^{k+1}(f(x))\,,\quad x\in \mathbb{R}^\star$$

   \item[(ii)] If $\alpha=k+r,\;  k\in \mathbb{N}\,,\;-\frac{1}{2}<r<\frac{1}{2},$
  $$ V_\alpha f
  (x)=c_1\,\int_{|x|}^{+\infty}
  (y^2-x^2)^{-r-\frac{1}{2}}\left(\frac{d}{dy^2}\right)^{k+1}(f)(sgn(x)y)\,y\,dy\,
  ,\;\,
  x\in\mathbb{R}^\star$$ where $\quad\displaystyle c_1=\frac{(-1)^{k+1}2\sqrt{\pi}}{\Gamma(\alpha+1)\Gamma(\frac{1}{2}-r)}$
\end{description}
\end{prop}
{\bf Proof.} Let $g\in \mathcal{E}(\mathbb{R})$ then we have\\
$\displaystyle <(\chi_\alpha^\star)^{-1}\,
T_f,g>\,=\,<(\chi_\alpha^{-1})^\star \,T_f,g>\, =\,
<T_f,\chi_\alpha^{-1}g>$ \\ {\bf First case $\displaystyle
\alpha=k+\frac{1}{2}\,,\;k\in\mathbb{N}$}\\ Invoking (ii) of
theorem \ref{qui-1}, we can write
\\ $\displaystyle <(\chi_\alpha^\star)^{-1}\,
T_f,g>\,=\,\frac{2^{2k+1}\,k!}{(2k+1)!}\int_\mathbb{R} f(x)\, x \,
\left(\frac{d}{dx^2}\right)^{k+1}(x^{2k+1}\, g(x))\,dx $\\
\hspace*{3cm}$\displaystyle=\,\frac{2^{2k+1}\,k!}{(2k+1)!}\,(I_1\,+\,I_2)$
 \\ where $$I_1=\int_0^\infty f(x)\, x \,
\left(\frac{d}{dx^2}\right)^{k+1}(x^{2k+1}\, g(x))\,dx$$ and
$$I_2=\int_{-\infty}^0f(x)\, x \,
\left(\frac{d}{dx^2}\right)^{k+1}(x^{2k+1}\, g(x))\,dx$$ By
integration by parts we have
\\ $\displaystyle I_1=\left[\frac{1}{2}f(x)
(\frac{d}{dx^2})^k(x^{2k+1}g)(x) \right]_0^{+\infty}\,-$
\\\hspace*{3,5cm}$\displaystyle\int_0^{+\infty}\left(\frac{d}{dx^2}\right)
f(x)\,\left(\frac{d}{dx^2}\right)^{k}(x^{2k+1}\, g(x))\,x\,dx$
\\ According to relation (\ref{a}) for $p=k$ and $m= 2k+1$, we
find that\\
$\displaystyle\left(\frac{d}{dx^2}\right)^{k}(x^{2k+1}\,g(x))=\sum_{i=0}^k
\beta_i^k x^{i+1}\,g^{(i)}(x) $.
\\ and $$I_1=-\int_0^{+\infty}\left(\frac{d}{dx^2}\right)
f(x)\,\left(\frac{d}{dx^2}\right)^{k}(x^{2k+1}\, g(x))\,x\,dx $$
Let $1\leq p\leq k+1$. After $p-1$ integrations by parts, we get
\\ $\displaystyle
I_1=\left[(-1)^{p-1}(\frac{d}{dx^2})^{p-1}(f(x))(\frac{d}{dx^2})^{k+1-p}(x^{2k+1}g(x))\right]_0^{+\infty}$
\\  \hspace*{3cm}$\displaystyle +
(-1)^p\int_0^{+\infty}\left(\frac{d}{dx^2}\right)^{p}
f(x)\,\left(\frac{d}{dx^2}\right)^{k-p+1}(x^{2k+1}\, g(x))\,x\,dx$
\\ Using relation (\ref{a}) again, we find
$$(\frac{d}{dx^2})^{p-1}(f(x))(\frac{d}{dx^2})^{k+1-p}(x^{2k+1}g(x))=x\sum_{i=0}^{k+1-p}\beta_i^{k+1-p}x^{i}g^{(i)}(x)
\sum_{i=0}^{p-1}\beta_i^{p-1}x^{i}f^{(i)}(x)$$ Since $f$ be in
$\mathcal{K}_0$, we obtain $$\displaystyle
I_1=(-1)^p\int_0^{+\infty}\left(\frac{d}{dx^2}\right)^{p}
f(x)\,\left(\frac{d}{dx^2}\right)^{k-p+1}(x^{2k+1}\,
g(x))\,x\,dx$$For $p=k+1$ we find $$\displaystyle
I_1=(-1)^{k+1}\int_0^{+\infty}\left(\frac{d}{dx^2}\right)^{k+1}
f(x)\, g(x)\,x^{2k+2}\,dx$$ As the same we establish that
$$\displaystyle
I_2=(-1)^{k+1}\int_{-\infty}^0\left(\frac{d}{dx^2}\right)^{k+1}
f(x)\, g(x)\,x^{2k+2}\,dx$$    Consequently,
$$<(\chi_\alpha^\star)^{-1}\,
T_f,g>\,=\,\frac{2^{2k+1}\,k!}{(2k+1)!}
(-1)^{k+1}\int_\mathbb{R}\left(\frac{d}{dx^2}\right)^{k+1} f(x)\,
g(x)\,x^{2k+2}\,dx$$ Which proves the wanted result for
$\displaystyle \alpha= k+\frac{1}{2}$.
\\ {\bf Second case $\displaystyle
\alpha=k+r\,,\;k\in\mathbb{N}\,,\;\frac{-1}{2}<r<\frac{1}{2}$ } \\
By virtue of (i) of theorem \ref{qui-1} and a change of variable,
we can write $$ \chi_\alpha^{-1}
g(x)=\frac{2\sqrt{\pi}x}{\Gamma(\alpha+1)\Gamma(\frac{1}{2}-r)}\left(\frac{d}{dx^2}\right)^{k+1}(x^{2k+1}\,h(x))
$$ where
$$h(x)=\int_0^1(1-u^2)^{-r-\frac{1}{2}}g(xu)\,u^{2\alpha+1}\,du$$
It's clear that $h\in\mathcal{E}(\mathbb{R})$, we proceed in a
similar way as in the first case, we just replace the function $g$
by the function $h$ and we obtain $$<(\chi_\alpha^\star)^{-1}\,
T_f,g>\,=\,c_1\int_\mathbb{R}\left(\frac{d}{dx^2}\right)^{k+1}
f(x)\, h(x)\,x^{2k+2}\,dx$$ Next, by a change of variable
$$\displaystyle <(\chi_\alpha^\star)^{-1}\, T_f,g>=c_1
\int_\mathbb{R}x(\frac{d}{dx^2})^{k+1} f(x)\left(\int_0^x
(x^2-t^2)^{-r-\frac{1}{2}}g(t)|t|^{2\alpha+1}dt\right)dx $$
\hspace*{3cm} $\displaystyle
=\,\frac{2\sqrt{\pi}\,(-1)^{k+1}}{\Gamma(\alpha+1)\Gamma(\frac{1}{2}-r)}(J_1+J_2)
$ \\ where $$J_1=\int_0^{+\infty}x(\frac{d}{dx^2})^{k+1}
f(x)\,\left(\int_0^x
(x^2-t^2)^{-r-\frac{1}{2}}g(t)\,|t|^{2\alpha+1}\,dt\right)\,dx $$
and $$J_2=\int_{-\infty}^0x(\frac{d}{dx^2})^{k+1}
f(x)\,\left(\int_0^x
(x^2-t^2)^{-r-\frac{1}{2}}g(t)\,|t|^{2\alpha+1}\,dt\right)\,dx $$
Applying Fubini's theorem at $J_1$ and $J_2$, we obtain $$J_1=
\int_0^{+\infty}\left(\int_t^{+\infty}
(x^2-t^2)^{-r-\frac{1}{2}}x(\frac{d}{dx^2})^{k+1} f(x)\,dx\right)
\, g(t)\,|t|^{2\alpha+1}\,dt $$ and $$J_2=
\int_{-\infty}^0\left(\int_{-\infty}^t
(x^2-t^2)^{-r-\frac{1}{2}}\, x(\frac{d}{dx^2})^{k+1}
f(x)\,dx\right)\,g(t)\,|t|^{2\alpha+1}\,dt $$ making a change of
variable in $J_2$ and using Chasles relation , we get
\\ $\displaystyle <(\chi_\alpha^\star)^{-1}\, T_f,g>=$  $$c_1
\int_\mathbb{R}\left(\int_{|t|}^{+\infty}
(x^2-t^2)^{-r-\frac{1}{2}}x(\frac{d}{dx^2})^{k+1}
f(sgn(t)x)\,dx\right) \, g(t)\,|t|^{2\alpha+1}\,dt $$
 \hfill $\square$
\begin{rmk} From proposition \ref{v} we deduce that the operators $V_\alpha$ and $\chi_\alpha^{-1}$ are
related by the following relation
\begin{equation}\label{vqui}
  \int_{\mathbb{R}}V_\alpha f(x)\, g(x) A(x)\, dx\,=\,  \int_{\mathbb{R}} f(x)\,\chi_\alpha^{-1} g(x)\, dx
\end{equation}
for all $f\in \mathcal{K}_0$ and $ g\in \mathcal{E}(\mathbb{R})$
\end{rmk}
\begin{lem}\label{inverse} Let $f$ be in
$\mathcal{D}(\mathbb{R})$. We have $W_\alpha(f)\in \mathcal{K}_0 $
and $\quad V_\alpha(W_\alpha(f)) = f$
\end{lem}
{\bf Proof.} For $f\in\mathcal{D}(\mathbb{R})$, by lemma
\ref{suppr} we have $W_\alpha(f)\in \mathcal{K}_0 $. Using
relations (\ref{vqui}) and (\ref{wqui}), we obtain for all $g\in
\mathcal{E}(\mathbb{R})$ $$\int_\mathbb{R}V_\alpha(W_\alpha
f)(x)\, g(x) A(x)\, dx\,=\int_\mathbb{R} f(x)\, g(x) A(x)\, dx $$
Thus $$V_\alpha(W_\alpha (f))(x)\,A(x) =f(x)\,A(x)\qquad a.e
\;\;x\in \mathbb{R}$$ Since $f\,A$ and $V_\alpha\circ W_\alpha
(f)\,A$ are both continuous functions on $\mathbb{R}^\star$ we
have $V_\alpha\circ W_\alpha (f)(x)=f(x)$ for all $x$ in
$\mathbb{R}^\star $ therefore $V_\alpha\circ W_\alpha (f)(x)=f(x)$
for all $x$ in $\mathbb{R}$. \hfill$\square$
\\

For $\displaystyle \alpha=k+\frac{1}{2}\;,\;\;k\in \mathbb{N}$, we
denote by $\Delta_{a,k+\frac{1}{2}}(\mathbb{R})$ the subspace of
$\mathcal{K}_0$ of functions $f$ infinitely differentiable on
$\mathbb{R}^\ast$ with support included in $[-a,a]$ verifying the
following condition :
$$\displaystyle(\frac{d}{dx^2})^{k+1}f\;\mathrm{ can \;be
\;extended \;to\; a\; function\; belonging\;
to\;}\mathcal{D}(\mathbb{R}).$$
\\ This space is provided with the topology defined by the semi
norms $\rho_n$ where $$\rho_n(f)=\sup_{\substack{ 0\leq p\leq n \\
x\in
[-a,a]}}\left|\big((\frac{d}{dx^2})^{k+1}f\big)^{(p)}(x)\right|\qquad
\quad ,\; n\in \mathbb{N} $$ We consider , for $k\in \mathbb{N}$,
the space $$\Delta_{k+\frac{1}{2}}(\mathbb{R})=\bigcup_{a\geq
0}\Delta_{a,k+\frac{1}{2}}(\mathbb{R})$$ endowed with the
inductive limit topology.
\begin{lem}\label{dw2}
For all $f$ in $\mathcal{D}_{a}(\mathbb{R})$ we have
$$\forall\,x\in \mathbb{R}^\ast,\qquad [W_{\frac{1}{2}}f]'(x)=\,-
x\,f(x)$$ and $$\forall\,\alpha
> \frac{1}{2}\,,\;\forall\;x\in\mathbb{R}^\ast,\qquad [W_\alpha
f]'(x)=-2\,\alpha x W_{\alpha-1}f(x)$$
\end{lem}
{\bf Proof.} Let $f\in \mathcal{D}_a(\mathbb{R})$\\ we have
$\displaystyle
W_{\frac{1}{2}}f(y)=\left\{\begin{array}{ll}\displaystyle
\frac{2\Gamma(\frac{3}{2})}{\sqrt{\pi}\Gamma(1)}\int_y^a
x\,f(x)\,dx & if\, y>0 \\
\displaystyle\frac{2\Gamma(\frac{3}{2})}{\sqrt{\pi}\Gamma(1)}
\int_{-y}^a x\,f(-x)\,dx & if\, y<0
\end{array}\right.$\\ We deduce that $\forall \,y\in \mathbb{R}^\ast\;,\;\;\displaystyle
[W_{\frac{1}{2}}f]'(y)=\,-y\,f(y)$

Now, we take $\displaystyle\alpha>\frac{1}{2}$, by lemma
\ref{suppr} $supp(W_\alpha f)\subset [-a,a]$. \\ Let $\varphi\in
\mathcal{D}(]0,+\infty[)$ then we have \\ \\ \hspace*{2cm}
$\displaystyle <[W_\alpha f]',\varphi>\,=\;-\,<W_\alpha
f,\varphi'>$\\ \hspace*{4,6cm}$\displaystyle\;
=-a_\alpha\int_0^a\int_y^a(x^2-y^2)^{\alpha-\frac{1}{2}}x\,f(x)\,dx\;\varphi'(y)\,dy$\\
%\hspace*{4,8cm}$\displaystyle =\lim_{\epsilon\rightarrow 0
%}-a_\alpha\int_\epsilon^a\int_y^a(x^2-y^2)^{\alpha-\frac{1}{2}}x\,f(x)\,dx\;\varphi'(y)\,dy$
\\ Using Fubini's theorem and an integration by parts we obtain %\\ $\displaystyle <[W_\alpha
%f]',\varphi>=\lim_{\epsilon\rightarrow 0
%}-a_\alpha\int_\epsilon^a\int_\epsilon^x(x^2-y^2)^{\alpha-\frac{1}{2}}\varphi'(y)\,dy\;x\,f(x)\,dx$
\\ $\displaystyle <[W_\alpha
f]',\varphi>=-a_\alpha\int_0^a\int_0^x(2\alpha-1)y(x^2-y^2)^{\alpha-\frac{3}{2}}\varphi(y)\,dy\;x\,f(x)\,dx$\\
Applying Fubini's theorem again we have \\ $\displaystyle
<[W_\alpha
f]',\varphi>=-2\alpha\int_0^ay\,W_{\alpha-1}f(y)\,\varphi(y)\,dy=<-2y\,W_{\alpha-1}
f,\varphi> $ \\ This proves that the derivative of the
distribution $W_\alpha f$ is the distribution defined by the
function $-2\alpha x\,W_{\alpha-1}$ on $]0,+\infty[$. The theorem
III p.54 in \cite{ls} allows us to say that the derivative on
$]0,+\infty[$ of the function $W_\alpha f$ is the function
$-2\alpha \,x W_{\alpha-1}f$.}\\ In the same way we obtain that
the derivative on $]-\infty,0[$ of the function $W_\alpha f$ is
the function $-2\alpha \,x W_{\alpha-1}f$. \hfill$\square$
\begin{thm}\label{th1} The operator $W_{k+\frac{1}{2}}$ is a topological isomorphism from
$\mathcal{D}_{a}(\mathbb{R})$ into
$\Delta_{a,k+\frac{1}{2}}(\mathbb{R})$ and its inverse is
$V_{k+\frac{1}{2}|\Delta_{a,k+\frac{1}{2}}(\mathbb{R})}$.
\end{thm}
{\bf Proof.} We will proceed by induction. According to lemma
\ref{dw2} we have
$W_{\frac{1}{2}}(\mathcal{D}_{a}(\mathbb{R}))\subset\Delta_{a,\frac{1}{2}}(\mathbb{R})$.
\\ Suppose for $k\in \mathbb{N}^\ast$ that
$W_{k-\frac{1}{2}}(\mathcal{D}_a(\mathbb{R})) \subset
\Delta_{a,k-\frac{1}{2}}(\mathbb{R})$. Let $f$ be in
$\mathcal{D}_a(\mathbb{R})$ then according to lemma \ref{dw2} and
the induction hypothesis, we deduce that
$\displaystyle(\frac{d}{dx^2})^{k+1}W_{k+\frac{1}{2}} f\in
\mathcal{D}(\mathbb{R})$. Furthermore from lemma \ref{suppr} we
conclude that $W_{k+\frac{1}{2}} f\in
\Delta_{a,k+\frac{1}{2}}(\mathbb{R})$.\\ In the other hand, let
$g\in \Delta_{a,k+\frac{1}{2}}(\mathbb{R})$. From proposition
\ref{v}, $\displaystyle V_{k+\frac{1}{2}}(g)$ can be extended to a
function in $\mathcal{D}_a(\mathbb{R})$. Since
$W_{k+\frac{1}{2}}(V_{k+\frac{1}{2}})(g)=g$ and by lemma
\ref{inverse} we deduce that $W_{k+\frac{1}{2}}$ is an isomorphism
from $\mathcal{D}_{a}(\mathbb{R})$ into
$\mathcal{D}_{a,k+\frac{1}{2}}(\mathbb{R})$ and its inverse is
$V_{k+\frac{1}{2}|\Delta_{a,k+\frac{1}{2}}(\mathbb{R})}$.\\ We
have $$\rho_n(W_{k+\frac{1}{2}}f)=\sup_{\substack{ 0\leq p\leq n
\\ x\in
[-a,a]}}\left|\big((\frac{d}{dx^2})^{k+1}W_{k+\frac{1}{2}}f\big)^{(p)}(x)\right|$$
From proposition \ref{v} and lemma \ref{inverse}
$$\rho_n(W_{k+\frac{1}{2}}f)=C\sup_{\substack{ 0\leq p\leq n \\
x\in [-a,a]}}|f^{(p)}(x)|$$ which proves the wanted result.
\hfill$\square$
\\

For $k\in \mathbb{N}$ we take $\alpha=k+r\,,\;r\in
]\frac{-1}{2},\frac{1}{2}[$.

We denote by $\Delta_{a,k+r}(\mathbb{R})$ the subspace of
$\mathcal{K}_0$ of functions $f$ infinitely differentiable on
$\mathbb{R}^\ast$ with support included in $[-a,a]$ verifying the
following condition :\\
$\displaystyle(\frac{d}{dx^2})^{k+1}\left(\int_1^{+\infty}(t^2-1)^{-r-\frac{1}{2}}f(xt)\,t\,dt\right)$
can be extended to a function belonging to
$|x|^{2r-1}\,\mathcal{D}(\mathbb{R})$

This space is provided with the topology defined by the semi norms
$q_n$ where $$q_n(f)=\sup_{\substack{ 0\leq p\leq n \\ x\in
[-a,a]}}\left|D^p\left(|x|^{-2r+1}(\frac{d}{dx^2})^{k+1}
\left(\int_1^{+\infty}(t^2-1)^{-r-\frac{1}{2}}f(xt)\,t\,dt\right)\right)\right|$$
We consider, for $k\in \mathbb{N}$, the space
$$\Delta_{k+r}(\mathbb{R})=\bigcup_{a\geq
0}\Delta_{a,k+r}(\mathbb{R})$$ endowed with the inductive limit
topology.
\begin{lem}\label{inverse3}
We have for all $f$ in $\Delta_{k+r}(\mathbb{R})$, $$
V_{k+r}(f)\in \mathcal{D}(\mathbb{R})\quad \mathrm{and} \quad
W_{k+r}(V_{k+r}(f))=f$$
\end{lem}
{\bf Proof.} Let $f\in \Delta_{k+r}(\mathbb{R})$,
\\ $\displaystyle
V_{k+r}(f)(y)=\int_{|y|}^{+\infty}(x^2-y^2)^{-r-\frac{1}{2}}\sum_{i=0}^{k+1}\beta_i^{k+1}x^i(sgn(y))^if^{(i)}(sgn(y)x)x\,dx$
\\ $\displaystyle
=|y|^{-2r+1}\int_1^{+\infty}(t^2-1)^{-r-\frac{1}{2}}\sum_{i=0}^{k+1}\beta_i^{k+1}y^it^if^{(i)}(yt)t\,dt$
\\ $\displaystyle
=|y|^{-2r+1}\sum_{i=0}^{k+1}\beta_i^{k+1}y^i(\frac{d}{dy})^{i}\left(\int_1^{+\infty}(t^2-1)^{-r-\frac{1}{2}}f{(yt)}\,t\,dt\right)$
\\ $\displaystyle
=|y|^{-2r+1}(\frac{d}{dy^2})^{k+1}\left(\int_1^{+\infty}(t^2-1)^{-r-\frac{1}{2}}f{(yt)}\,t\,dt\right)$
\\ Then $V_{k+r}(f)\in\mathcal{D}(\mathbb{R})$ and from relation
(\ref{wqui}), for all $g\in \mathcal{E}(\mathbb{R})$, \\
$$\displaystyle
\int_{\mathbb{R}}W_{k+r}(V_{k+r}(f))(x)\,g(x)\,dx=\int_{\mathbb{R}}V_{k+r}(f)(x)\,\chi_{k+r}(g(x))A(x)\,dx$$
By relation(\ref{vqui}), we deduce $$\displaystyle
\int_{\mathbb{R}}W_{k+r}(V_{k+r}(f))(x)\,g(x)\,dx=\int_{\mathbb{R}}f(x)\,g(x)\,dx$$
\\ Therefore $$W_{k+r}(V_{k+r}(f))(x)=f(x)\;,\qquad a.e.\;x\in
\mathbb{R}$$ Since $W_{k+r}(V_{k+r}(f))$ and $f$ are both
continuous functions on $\mathbb{R}^\ast$, we get $\forall\;x\in
\mathbb{R}^\ast\;,\qquad W_{k+r}(V_{k+r}(f))(x)=f(x)$
 \hfill$\square$
\begin{thm}\label{th2}
$W_{k+r}$ is a topological isomorphism from
$\mathcal{D}_{a}(\mathbb{R})$ into $\Delta_{a,k+r}(\mathbb{R})$
and its inverse is $V_{k+r|\Delta_{a,k+r}(\mathbb{R})}$
\end{thm}
{\bf Proof.} Let $f\in \mathcal{D}_a(\mathbb{R})$, from
proposition \ref{suppr} and lemma \ref{inverse}, we can write
$$|x|^{-2r+1}\displaystyle(\frac{d}{dx^2})^{k+1}\left(\int_1^{+\infty}(t^2-1)^{-r-\frac{1}{2}}
W_{k+r}(f)(xt)\,t\,dt\right) =V_{k+r}(W_{k+r}(f))=f$$ which proves
that $W_{k+r}(f) \in \Delta_{a,k+r}(\mathbb{R})$.\\ Furthermore,
by lemma \ref{inverse3} and lemma \ref{inverse}, one can deduce
that $W_{k+r}$ is bijective and
$V_{k+r|\Delta_{a,{k+r}}(\mathbb{R})}$ is its inverse on
$\Delta_{a,{k+r}}(\mathbb{R})$

Now, the fact that $$q_n(W_{k+r} (f))=Cp_n(f)$$ allows us to
conclude that $W_{{k+r}}$ is a topological isomorphism from
$\mathcal{D}_{a}(\mathbb{R})$ into $\Delta_{a,{k+r}}(\mathbb{R})$.
\hfill$\square$ \\ The following theorem is a consequence of
theorem \ref{th1} and theorem \ref{th2} .
 \begin{thm}
$W_{\alpha}$ is a topological isomorphism from
$\mathcal{D}(\mathbb{R})$ into $\Delta_{\alpha}(\mathbb{R})$ and
its inverse is $V_{\alpha|\Delta_{\alpha}(\mathbb{R})}$
\end{thm}
\section{Paley Wiener type theorem associated to Bessel-Struve transform }
In this section we will try to characterize the range of
$\mathcal{D}(\mathbb{R})$ by Bessel-Struve transform.
\subsection{Some properties of Bessel-Struve transform on $\mathcal{D}(\mathbb{R})$}
\begin{prop}\label{pw}
Let $f$ be a function  in $\mathcal{K}_0$. Then the classical
Fourier transform $\mathcal{F}(f)$ can be extended to an analytic
function on $\mathbb{C}$ that we denote again $\mathcal{F}(f)$ and
we have
\begin{equation}\label{zkf}
\forall k\in \mathbb{N}^\star\,,\qquad
(iz)^k\mathcal{F}(f)(z)=\int_\mathbb{R}f^{(k)}(x)(e^{-ixz}-\sum_{n=0}^{k-1}\frac{(-izx)^n}{n!})\,dx
\end{equation}
\end{prop}
{\bf Proof.} Let $f\in \mathcal{K}_0$ and $x\in \mathbb{R}$. We
have, for all $z\in \mathbb{C}$ such that $Im(z)<b$,
$$|f(x)e^{-ixz}|\leq |f(x)|e^{|x|b}$$ Since $f$ has a bounded
support and the function $:z\longrightarrow e^{-ixz}\,f(x)$ is
analytic on $\mathbb{C}$, we conclude that $\mathcal{F}(f)$ is
analytic function on $\mathbb{C}$.\\ To prove relation
(\ref{zkf}), we proceed by induction. \\ For $k=1$ we can write $$
iz \mathcal{F}(f)(z)= \lim_{\epsilon\rightarrow
0}\left(\int_{-\infty}^{-\epsilon}f(x)ize^{-ixz}dx+\int_{\epsilon}^{+\infty}f(x)ize^{-ixz}dx\right)$$
By integration by parts, we obtain $$ iz \mathcal{F}(f)(z)=
\lim_{\epsilon\rightarrow
0}\left(\int_{-\infty}^{-\epsilon}f'(x)(e^{-ixz}-1)dx+\int_{\epsilon}^{+\infty}f'(x)(e^{-ixz}-1)dx\right)$$
Therefore $$iz
\mathcal{F}(f)(z)=\int_\mathbb{R}f'(x)(e^{-ixz}-1)dx$$ We suppose
that we have, for $k\in\mathbb{ N}^\ast$ $$
(iz)^k\mathcal{F}(f)(z)=\int_\mathbb{R}f^{(k)}(x)(e^{-ixz}-\sum_{n=0}^{k-1}\frac{(-izx)^n}{n!})\,dx$$
Then \\
$\displaystyle(iz)^{k+1}\mathcal{F}(f)(z)=iz(iz)^k\mathcal{F}(f)(z)$\\
\hspace*{2,65cm} $=\displaystyle \lim_{\epsilon\rightarrow
0}\left(\int_{-\infty}^{-\epsilon}f^{(k)}(x)iz(e^{-ixz}-\sum_{n=0}^{k-1}\frac{(-izx)^n}{n!})\,dx\right.
$\\ \hspace*{5cm} $\left.\displaystyle +
\int_{\epsilon}^{+\infty}f^{(k)}(x)iz(e^{-ixz}-\sum_{n=0}^{k-1}\frac{(-izx)^n}{n!})\,dx
\right)$\\ By integration by parts, we find $$
(iz)^{k+1}\mathcal{F}(f)(z)=\int_\mathbb{R}f^{(k+1)}(x)(e^{-ixz}-\sum_{n=0}^{k}\frac{(-izx)^n}{n!})\,dx$$
which completes the proof. \hfill$\square$
\begin{rmk}
we have for $f$ in $\mathcal{K}_0$, $\displaystyle
\int_\mathbb{R}f^{(k)}(t)e^{-itz}dt$ is generally divergent. The
relation (\ref{zkf}) means that we have for $f$ in $\mathcal{K}_0$
$(iz)^k\mathcal{F}(f)$ is the finite part of the divergent
integral citing below. (For the definition of the finite part of
the divergent integral, one can see \cite{ls}).
\end{rmk}
\begin{thm}
Let $a>0$ and $f$ a function in $\mathcal{D}_a(\mathbb{R})$ then
$\mathcal{F}_{B,S}^\alpha(f)$ can be extended to an analytic
function on $\mathbb{C}$ that we denote again
$\mathcal{F}_{B,S}^\alpha(f)$ verifying
\begin{equation}\label{01}
\forall k\in \mathbb{N}^\star\,,\qquad
|\mathcal{F}_{B,S}^\alpha(f)(z)|\leq C e^{a|z|}
\end{equation}
\end{thm}
{\bf Proof.} Using proposition \ref{pw}, corollary \ref{xnw} and
relation (\ref{wf}), we find that $\mathcal{F}_{B,S}^\alpha(f)$
can be extended to an analytic function on $\mathbb{C}$ and
\begin{equation}\label{11}
\forall k\in \mathbb{N}^\star\,,\qquad
(iz)^k\mathcal{F}_{B,S}^\alpha(f)(z)=\int_\mathbb{R}(W_\alpha
f)^{(k)}(x)\sum_{n=k}^{+\infty}\frac{(-izx)^n}{n!}\,dx
\end{equation} Then $$ |z|^k|\mathcal{F}_{B,S}^\alpha(f)(z)|\leq
|z|^k \int_\mathbb{R}|x^k\,(W_\alpha
f)^{(k)}(x)|\;\sum_{n=0}^{+\infty}\frac{|zx|^n}{n!}\,dx $$
Therefore $$\mathrm{if} \; z\neq0\,,\quad
|\mathcal{F}_{B,S}^\alpha(f)(z)|\leq \int_{-a}^a \,|x^k\,(W_\alpha
f)^{(k)}(x)|\,dx \;.\;e^{a|z|}$$ and $$\mathrm{if} \; z=0\,,\quad
|\mathcal{F}_{B,S}^\alpha(f)(0)|\leq \|f\|_{1,\alpha}$$ So, we
obtain the relation (\ref{01}).\hfill$\square$
\begin{prop}
Let $a>0$ and $f$ a function in $\mathcal{D}_a(\mathbb{R})$ then
we have
\begin{equation}\label{02}
\forall n\in\mathbb{N}\,,\quad
[\mathcal{F}_{B,S}^\alpha(f)]^{(n)}\in C_0(\mathbb{R})
\end{equation}
\end{prop}
{\bf Proof.}  Using The relation (\ref{wf}) and the derivation
theorem we obtain the following equality
$$[\mathcal{F}_{B,S}^\alpha(f)]^{(n)}=\mathcal{F}((-it)^nW_\alpha
f)$$ Since $W_\alpha f\in \mathcal{K}_0$ then $t^nW_\alpha f$
belongs to $L^1(\mathbb{R})$ and the relation (\ref{02}) can be
deduced.\hfill$\square$
 %%%%%%%%%%%%%%%%%%%%%%%%%%%%%%%%%%%%%%%%%%%%%%%%%%%%%%%%%%%%%%%%%%
\subsection{Range of $\mathcal{D}(\mathbb{R} )$ by Bessel-Struve transform
for half integers}
 Let $a>0$, $\mathcal{H}_a$ designates the space
of entire functions $f$ verifying : $$\forall n\in
\mathbb{N}\,,\;\exists\,c_n>0\,;\;\forall z\in \mathbb{C}\,,\qquad
(1+|z|^2)^n |f(z)|e^{-a\, Im(z)}\,<\,c_n$$ and
$$\mathcal{H}=\bigcup_{a>0}\mathcal{H}_a$$ We introduce the space
$\Lambda_{a,\frac{1}{2}}$ the space of entire functions $g$
verifying
  \begin{equation}\label{cond} \exists\,h\in \mathcal{H}_a\;\forall \; z\in
  \mathbb{C}^\ast \qquad g(z)=\frac{h'(z)-h'(0)}{z}
  \end{equation}
and we denote $\displaystyle
\Lambda_{\frac{1}{2}}=\bigcup_{a>0}\Lambda_{a,\frac{1}{2}}$
\begin{thm}\label{pdemi} We have
$$\mathcal{F}^{\frac{1}{2}}_{BS}(\mathcal{D}(\mathbb{R}))=\Lambda_{\frac{1}{2}}$$
\end{thm}
{\bf Proof.} Let $f\in \mathcal{D}(\mathbb{R})$ and $z\in
\mathbb{C}$. Using relation (\ref{11}), we have
$$\mathcal{F}^{\frac{1}{2}}_{BS}(f)(z)=\lim_{\substack{\epsilon\rightarrow
0 \\ \epsilon
>0}}\left(\int_{-\infty}^{-\epsilon}W_{\frac{1}{2}}f(t)e^{-izt}dt+
\int_{\epsilon}^{+\infty}W_{\frac{1}{2}}f(t)e^{-izt}dt\right)$$ By
integration by parts
$$-iz\,\mathcal{F}^{\frac{1}{2}}_{BS}(f)(z)=-c+\mathcal{F}(xf)(z)$$
where $$c=\lim_{\substack{x\rightarrow 0 \\
x>0}}W_{\frac{1}{2}}f(x)-\lim_{\substack{x\rightarrow 0 \\
x<0}}W_{\frac{1}{2}}f(x)$$ Furthermore
$\mathcal{F}(x\,f)(z)=i[\mathcal{F}(f)]'(z)$ and for $z=0$, we get
$c=i[\mathcal{F}(f)]'(0)$. Since $f\in\mathcal{D}(\mathbb{R}),\;\;
\mathcal{F}(f)\in \mathcal{H}$ and $$\displaystyle
\mathcal{F}^{\frac{1}{2}}_{BS}(f)(z)=\frac{[\mathcal{F}(-f)]'(z)-[\mathcal{F}(-f)]'(0)}{z}$$
Now let $g$ an entire function verifying relation (\ref{cond}).\\
From classical Paley-Wiener theorem we have $$\exists \,f\in
\mathcal{D}_a(\mathbb{R}) \;\mathrm{\;such\: that\;}
h(z)=\mathcal{F}(f)(z)$$ Then
$$[\mathcal{F}(f)]'(z)-[\mathcal{F}(f)]'(0)=z\,g(z)$$ Therefore,
for $\lambda \neq0$\\ $\displaystyle
g(\lambda)\;=\,-\frac{i}{\lambda}(\mathcal{F}(t\,f)(\lambda)-\mathcal{F}(t\,f)(0))$\\
\hspace*{1cm}$\displaystyle=\int_\mathbb{R}
f(t)(\frac{-\sin(\lambda\, t)}{\lambda\,t}+i\frac{1-\cos(\lambda\,
t)}{\lambda\, t})t^2\,dt$\\
\hspace*{1cm}$\displaystyle=-\,\int_\mathbb{R}
f(t)S_{\frac{1}{2}}(-i\,\lambda\, t)\,t^2\,dt$\\
\hspace*{1cm}$\displaystyle=\mathcal{F}^\frac{1}{2}_{BS}(-f)(\lambda)$
 \hfill$\square$

By induction, we can build the range of $\mathcal{D}(\mathbb{R})$
by $\mathcal{F}_{B,S}^{k+\frac{1}{2}}$ from theorem \ref{pdemi}
and the following proposition.
\begin{prop}\label{pfalfa}
For $\displaystyle\alpha>\frac{1}{2}$, the following assertions
are equivalent
\begin{description}
  \item[(i)] $g=\mathcal{F}_{BS}^{\alpha}(f)\, \mathrm{where} \;f\in \mathcal{D}_a(\mathbb{R})$
  \item[(ii)] $g$ is extented to an entire function $\tilde{g}$
  verifying
  \begin{equation}\label{cond2} \exists \; h\in
  \mathcal{F}_{BS}^{\alpha-1}(\mathcal{D}_a(\mathbb{R}))\;;\;\;\forall \, z\in
  \mathbb{C}\quad \quad \tilde{g}(z)=\alpha\,\frac{h'(z)-h'(0)}{z}
  \end{equation}
\end{description}
\end{prop}
{\bf Proof.}  Let $f\in \mathcal{D}(\mathbb{R})$ and $z\in
\mathbb{C}$. We proceed in a similar way as in theorem \ref{pdemi}
and we obtain $$iz
\mathcal{F}^{\alpha}_{BS}(f)(z)-c=\mathcal{F}([W_{\alpha}f]')(z)$$
where $$c=\lim_{\substack{x\rightarrow 0 \\
x>0}}W_{\frac{1}{2}}f(x)-\lim_{\substack{x\rightarrow 0 \\
x<0}}W_{\frac{1}{2}}f(x)$$ Furthermore, \;
$\displaystyle\mathcal{F}([W_{\alpha}f]')(x)=-i\,\alpha[\mathcal{F}(W_{\alpha-1}f)]'(x)$.
Since $f\in\mathcal{D}_a(\mathbb{R})\,,$\\
$\mathcal{F}(W_{\alpha-1}f)\in
\mathcal{F}^{\alpha-1}_{BS}(\mathcal{D}(\mathbb{R}))$ and
$\mathcal{F}^{\alpha}_{BS}(f)$ verifies relation (\ref{cond2}).\\
Now let $g$ be an entire function verifying relation
(\ref{cond2}).\\ Then $$\exists \,f\in \mathcal{D}_a(\mathbb{R})
\;\;such\: that\;
z\,g(z)-c=\alpha[\mathcal{F}^{\alpha-1}_{BS}(f)]'(z)$$ Therefore
$$\alpha([\mathcal{F}^{\alpha-1}_{BS}(f)]'(z)-[\mathcal{F}^{\alpha-1}_{BS}(f)]'(0))=z\,g(z)$$
In the other hand
$\displaystyle\mathcal{F}^{\alpha}_{BS}(f)(z)=\frac{1}{iz}\big(\lim_{\substack{x\rightarrow
0 \\ x>0}}W_{\frac{1}{2}}f(x)-\lim_{\substack{x\rightarrow 0 \\
x<0}}W_{\frac{1}{2}}f(x)-\alpha
[\mathcal{F}^{\alpha-1}_{BS}f]'\big)$\\ then we take $z=0$ we have
$$-\alpha\,[\mathcal{F}^{\alpha-1}_{BS}(f)]'(0)=\lim_{\substack{x\rightarrow
0 \\ x>0}}W_{\frac{1}{2}}f(x)-\lim_{\substack{x\rightarrow 0 \\
x<0}}W_{\frac{1}{2}}f(x)$$ Finally, by integration by parts we get
$$\displaystyle
g(z)=\frac{1}{z}\left((\lim_{\substack{x\rightarrow 0 \\
x>0}}W_{\frac{1}{2}}f(x)-\lim_{\substack{x\rightarrow 0 \\
x<0}}W_{\frac{1}{2}}f(x))-\alpha
[\mathcal{F}^{\alpha-1}_{BS}f]'\right)=i\mathcal{F}^{\alpha}_{BS}(f)(z)$$
 \hfill$\square$
\subsection{Schwartz Paley Wiener theorem}
In this subsection we will prove a Paley Wiener theorem in
distributions space with bounded support. (we can cite these
references : \cite{ls} \cite{vk2})
\begin{df}
We define the Fourier Bessel-Struve Transform on
$\mathcal{E}'(\mathbb{R})$ by \begin{equation}\forall \, T\in
\mathcal{E}'(\mathbb{R})\,,\; \mathcal{F}
_{B,S}^\alpha(T)(\lambda)\,=\,<T,S(-i\lambda.)>
\end{equation}
\end{df}
\begin{prop}
For all $T\in \mathcal{E}'(\mathbb{R})\,,$
\begin{equation}
\mathcal{F}_{B,S}^\alpha(T)\,=\;\mathcal{F}\circ
\chi_\alpha^\star(T)
\end{equation}
\end{prop}
{\bf Proof.} We have
$\mathcal{F}_{B,S}^\alpha(T)(\lambda)\,=\,<T,S(-i\lambda.)>\,=\,<T,\chi_\alpha(exp(-i\lambda.))>$
\hspace*{5,5cm}$=\,<\chi_\alpha^\star
T,exp(-i\lambda.)>\,=\;\mathcal{F}\circ \chi_\alpha^\star(T)$
  \hfill$\square$
\begin{lem}\label{tqui1}
Let $T\in \mathcal{E}'(\mathbb{R})$ , then $$supp(T)\subseteq
[-b,b]\Longleftrightarrow supp(\chi_\alpha^\star(T))\subseteq
[-b,b]$$
\end{lem}
{\bf Proof.} Let $T\in \mathcal{E}'(\mathbb{R})$ such that
$supp(T)$ included in $[-b,b]$ . For
$\varphi\in\mathcal{D}(\mathbb{R})$ with support in $[-b,b]^c$ we
have $\chi_\alpha \varphi$ have the support included in $[-b,b]^c$
therefore $<\chi_\alpha^\star
(T),\varphi>=<T,\chi_\alpha\varphi>=0$ then $\chi_\alpha^\star
(T)$ have the support in $[-b,b]$
\\ Now we consider a distribution $T$ such that
$supp\,(\chi_\alpha^\star(T))$ included in $[-b,b]$. For
$\varphi\in\mathcal{D}(\mathbb{R})$ with support in $[-b,b]^c$ we
have $$<T,\varphi>=<(\chi_\alpha^{-1})^\star\circ\chi_\alpha^\star
(T),\varphi>=<\chi_\alpha^\star (T), \chi_\alpha^{-1}\varphi>$$
Using theorem \ref{qui-1} $supp\,(\chi_\alpha^{-1}\varphi)$
included in $[-b,b]$ so $<T,\varphi>=0$ which complete the proof.
\hfill$\square$
\begin{thm} Let $b>0$ and $f\in \mathcal{E}(\mathbb{R})$.There is an equivalence between the two following assertions
\begin{enumerate}
  \item There exists a distribution
$T\in\mathcal{E}'(\mathbb{R})$ with support included in $[-b,b]$
such that $f\,=\,\mathcal{F}_{B,S}(T)$
  \item $f$ is extended to an
analytic function $\tilde{f}$ on $\mathbb{C}$ such
  that
  \begin{equation}\label{2}
  \exists m \in \mathbb{N}\,,\;\exists  c>0\,,\; \forall z\in \mathbb{C}\;\;
  |\tilde{f}(z)|\,\leq\, c\, (1+|z|^2)^{\frac{m}{2}}\,e^{b(Im(z))}
  \end{equation}
\end{enumerate}
\end{thm}
{\bf Proof.} Let $T\in \mathcal{E}'(\mathbb{R})$ such that
$supp(T)$ included in $[-b,b]$ then from lemma \ref{tqui1} and the
classical Paley-Wiener Schwartz , (one can see \cite{vk2}), on the
distribution $\chi_\alpha^\star(T)$ we obtain
$\mathcal{F}_{B,S}^\alpha $ verifies 2).
\\ We suppose that  $f$ is extended to an analytic function $\tilde{f}$ on
$\mathbb{C}$ verifying (\ref{2}) then there exists a distribution
$W$ with support included in $[-b,b]$ such that
$f=\mathcal{F}(W)$. Using corollary \ref{tqui-1} we can deduce the
wanted result.
 \hfill$\square$

\bibliographystyle{unsrt}

\begin{thebibliography}{11}
\bibitem{gs} A.Gasmi and M.Sifi, The Bessel-Struve Intertwining Operator on $\mathbb{C}$ and Mean-periodic Functions,
Hindawi Publishing Corp,IJMMS 2004:59,3171-3185.
\bibitem{gs2} A.Gasmi and M.Sifi, Analytic Mean-Periodic Functions Associated with the
Bessel-Struve Operator on a Disk,Global Journal of Pure and
Applied Mathematics. ISSN 0973-1768 Vol.1 No.1(2005), pp. 55-68.
\bibitem{kam} L. Kamoun and M. Sifi, Bessel-Struve Intertwining
Operator and Generalized Taylor Series on the Real Line, Integral
Transforms and Special Functions, Vol.16, January 2005, 39-55.
\bibitem{vk} Vo-Khac Kohan, Distributions Analyse de Fourier
Op\'erateurs aux D\'eriv\'ees Partielles, TomeI, Librairie
Vuiebert 1972
\bibitem{vk2} Vo-Khac Kohan, Distributions Analyse de Fourier
Op\'erateurs aux D\'eriv\'ees Partielles, TomeII, Librairie
Vuiebert 1972
\bibitem{ls} L.Schwartz, Th\'eorie des Distributions, Hermann Paris
1966
\bibitem{t1} K. Trim\`eche, Generalised Harmonic Analysis and Wavelets
Packets, Gordon and Breach Science Publishers, 2001
\bibitem{t2} K. Trim\`eche, Transformation Int\'egrale de Weyl et
Th\'eor\`eme de Paley Wiener Associ\'es \`a un Opérateur
Diff\'erentiel Singulier sur $(0,+\infty)$, J.Math.pures et appli.
60,(1981), 51-98
\bibitem{t3} K. Trim\`eche, Inversion Formulas for the Dunkl Intertwining Operator and
Its Dual on Spaces of Functions and Distributions, SIGMA 4 , 2008,
087, 22 p.
\bibitem{t4}  K. Trim\`eche, A New Harmonic Analysis Related To The Dunkl Operators Theory,
preprint
\bibitem{Wa} G.N. Watson, A Treatise on the Theory of Bessel Functions, Camb. Univ. Press,
Cambridge, 1966.
\end{thebibliography}

\end{document}